\documentclass[12pt]{amsart}
\usepackage[T1]{fontenc}
\usepackage[a4paper,margin=0.7in]{geometry}
\usepackage{amssymb}
\newtheorem{theorem}{Theorem}[section]
\newtheorem{lemma}[theorem]{Lemma}
\newtheorem{proposition}[theorem]{Proposition}
\newtheorem{corollary}[theorem]{Corollary}
\newtheorem{definition}[theorem]{Definition}
\newtheorem{fact}[theorem]{Fact}
\newtheorem{remark}[theorem]{Remark}
\newcommand{\N}{\mathbb{N}}
\newcommand{\dd}{\mathbf{d}}
\newcommand{\rev}[1]{\widetilde{#1}}
\newcommand{\pal}[1]{#1^{(+)}}
\usepackage[colorlinks,
linkcolor=blue,
anchorcolor=blue,
citecolor=blue
]{hyperref}

\author{J. Hamoud}
\address{\textbf{Jasem Hamoud:}  Moscow Institute of Physics and Technology}
\email{hamoud.math@gmail.com}
\thanks{}

\title[Lengths and Counting of Palindromic Prefixes]{Lengths and Counting of Palindromic Prefixes of Standard Sturmian Words}

\date{}
\begin{document}

\begin{abstract}
Every standard Sturmian word $c_\dd$ is determined by a directive sequence $\dd=(d_i)_{i\ge 1}$ of positive integers, and its palindromic prefixes are obtained by iterated palindromic closure. We describe the lengths of these prefixes by the numbers $j\,\ell_{i-1}+\ell_{i-2}-2$, where $\ell_n$ is the length of the $n$-th standard word and $1\le j\le d_i$, and we show that each nonzero length is attained by exactly one pair $(i,j)\neq(1,1)$. For a purely periodic directive sequence of period $k$ we then prove that the number $\pi_\dd(N)$ of nonempty palindromic prefixes of length at most $N$ satisfies
\[
\pi_\dd(N)\sim\frac{d_1+\cdots+d_k}{\log\lambda_\dd}\,\log N,
\]
where $\lambda_\dd$ is the largest eigenvalue of the matrix $M(d_1)\cdots M(d_k)$. For the Fibonacci word this gives $\pi_\dd(N)\sim\log N/\log\varphi$, where $\varphi=(1+\sqrt5)/2$.
\end{abstract}

\maketitle

\noindent\rule{\textwidth}{0.4pt}

\noindent\textbf{Keywords:} Sturmian words, standard Sturmian words, palindromic prefixes, palindromic closure, directive sequence, continued fractions, asymptotic counting.

\medskip

\noindent{\bf MSC 2020:} 68R15, 11A55, 37B10.

\noindent\rule{\textwidth}{0.4pt}

\section{Introduction}

Standard Sturmian words are generated by iterated palindromic closure, so their palindromic prefixes form the backbone of the word. 
They have recently appeared in several lines of work.
Dvo\v{r}\'ak\'ov\'a and Hendrychov\'a showed that the palindromic prefixes of standard Sturmian sequences have string attractors of size two, and used this to describe string attractors of pseudopalindromic prefixes of Rote sequences \cite{DH24}. 
Dvo\v{r}\'ak\'ov\'a then described string attractors of all factors of episturmian sequences, again through palindromic closures \cite{Dvo22}. 
Fici, Shallit and Simpson proved that every nonempty prefix of a standard episturmian word, and every nonempty factor of a Sturmian word, is a palindromic periodicity \cite{FSS24}. 
In a different direction, palindromic factorizations of prefixes and factors of aperiodic words have been studied through the palindromic length: Rukavicka investigated how powers can be reduced in an infinite word while keeping the palindromic length bounded \cite{Ruk22,Ruk24}.

The palindromic prefixes of standard Sturmian words arise from the palindromic closure $\pal{w}$ of a finite word $w$, that is, the shortest palindrome having $w$ as a prefix. 
Droubay, Justin and Pirillo \cite{DJP} used iterated palindromic closure along a directive word to define standard episturmian words over an arbitrary alphabet. 
Justin and Pirillo \cite{JP} showed that the palindromic prefixes of a standard episturmian word are exactly the words $u_n$, $n>0$, defined by $u_{n+1}=\pal{(u_nx_n)}$, where $x_1x_2\cdots$ is the directive word. 
They also showed that this recurrence can be written as $u_{n+1}=u_nx_nu_n$ if $x_n$ does not occur in $x_1\cdots x_{n-1}$, and $u_{n+1}=u_nu_{P(n)}^{-1}u_n$ otherwise, with $P(n)$ the last earlier position of $x_n$ \cite{JP}. 
On a two-letter alphabet the directive word $0^{d_1}1^{d_2}0^{d_3}\cdots$ gives the characteristic Sturmian word of slope $[0;d_1+1,d_2,d_3,\dots]$ \cite{Fischler}.

Previous studies have addressed various aspects of Sturmian words and their prefixes; however, since a comprehensive review lies beyond the scope of this work, the reader is referred to the literature for further details~\cite{Fischler2007,DeLuca2013, Pel2015,DeLuca2017, Gabric2022, Rampersad2025}.

In these works the palindromic prefixes are used through their structure. The present paper asks a quantitative question instead: how many palindromic prefixes of a standard Sturmian word have length at most $N$, and how does this number grow with $N$?

\medskip

The goal of the paper is to study the palindromic prefixes of standard Sturmian words $c_\dd$ through their directive sequences $\dd=(d_1,d_2,\dots)$. We first rederive the lengths of these prefixes as the numbers $j\,\ell_{i-1}+\ell_{i-2}-2$, indexed by pairs $(i,j)$ with $1\le j\le d_i$, where $\ell_n$ is the length of the $n$-th standard word. We then use this description to count them: for a purely periodic directive sequence $\dd=(d_1,\dots,d_k)^\infty$, we determine the asymptotic behaviour of the number $\pi_\dd(N)$ of nonempty palindromic prefixes of length at most $N$ in terms of $d_1+\cdots+d_k$ and the largest eigenvalue $\lambda_\dd$ of $M(d_1)\cdots M(d_k)$.

This paper is organized as follows. Section~\ref{sec2per} fixes the notation and recalls Sturmian and standard words, the construction of $c_\dd$ by iterated palindromic closure, and the length set of the palindromic prefixes. Section~\ref{mainsec03} contains the main results. Subsection~\ref{sec01A1} derives the length set from the closure construction and shows that each nonzero length corresponds to exactly one pair $(i,j)\ne(1,1)$. Subsection~\ref{sec00Asymptotic} proves the asymptotic formula for $\pi_\dd(N)$ in Theorem~\ref{Thm01Asymptotic}. The last two sections contain a discussion of the result and the conclusion.
\section{Preliminaries}\label{sec2per}
\paragraph{Notation and conventions.}
We work over the binary alphabet $A=\{0,1\}$. We denote by $A^*$ the set of finite words and by $\varepsilon$ the empty word. The length of a word $w$ is written $|w|$, and $|w|_a$ denotes the number of occurrences of the letter $a$ in $w$. A word $u$ is a \emph{factor} of $w$ if $w=xuy$ for some words $x,y$; it is a \emph{prefix} if $x=\varepsilon$ and a \emph{suffix} if $y=\varepsilon$. For a word $w$ and an integer $j\ge 0$, $w^j$ is the word $w$ repeated $j$ times, with $w^0=\varepsilon$. The \emph{reversal} of $w=w_1\cdots w_n$ is $\rev{w}=w_n\cdots w_1$, and $w$ is a \emph{palindrome} if $w=\rev{w}$, the empty word and single letters being palindromes. A prefix of a word $x$ that is a palindrome is called a \emph{palindromic prefix} of $x$.
Regarding asymptotics, $f(N)\sim g(N)$ means that $f(N)/g(N)\to 1$ as $N\to\infty$, and all logarithms are natural.

\subsection{Sturmian words and standard words}\label{secsubstu001}
\begin{definition}[\cite{Glen2006}]\label{defstu0001}
The complexity of an infinite word $x$ is $p_x(n)$, the number of distinct factors of length $n$. An infinite word is aperiodic if it is not ultimately periodic.
\end{definition}

\begin{fact}[\cite{Lothaire}]\label{factstu0001}
An aperiodic word $x$ satisfies $p_x(n)\ge n+1$ for all $n$.
\end{fact}

\begin{definition}[\cite{DeLuca}]\label{defstu0002}
An infinite word $x$ is Sturmian if $p_x(n)=n+1$ for all $n\ge 0$.
\end{definition}

\begin{fact}[\cite{DeLuca}]\label{factstu0002}
An infinite word $x$ is Sturmian if and only if it is aperiodic and balanced, that is, $\bigl||u|_1-|v|_1\bigr|\le 1$ for any two factors $u,v$ of $x$ of equal length.
\end{fact}

\begin{fact}[\cite{Lothaire}]\label{factstu0003}
Every Sturmian word has a frequency $\alpha$ of the letter $1$, which is irrational and is called its slope. The set of finite factors of a Sturmian word depends only on $\alpha$.
\end{fact}

\begin{definition}[\cite{Lothaire}]\label{defstu0004}
Given an irrational $\alpha\in(0,1)$, the characteristic word $c_\alpha$, called the standard word of slope $\alpha$, is the unique Sturmian word of slope $\alpha$ all of whose prefixes $u$ have both $0u$ and $1u$ as factors.
\end{definition}

The paper uses ``standard Sturmian word'' for $c_\alpha$.

\begin{definition}\label{defstu0005}
A directive sequence is $\dd=(d_1,d_2,d_3,\dots)$ with $d_1\ge 0$ and $d_i\ge 1$ for $i\ge 2$. The standard words are
\[
s_{-1}=1,\qquad s_0=0,\qquad s_n=s_{n-1}^{\,d_n}\,s_{n-2}\quad(n\ge 1).
\]
The word $s_n$ is a prefix of $s_{n+1}$ for $n\ge 1$.
\end{definition}
From now on $d_1\ge 1$, so that $d_i\ge 1$ for all $i\ge 1$.
The sequence $(s_n)$ converges to an infinite word $c_\dd$, which is a standard Sturmian word, and every standard Sturmian word is obtained in this way \cite{DeLuca}. 
The lengths are denoted $\ell_n=|s_n|$, with $\ell_{-1}=\ell_0=1$ and $\ell_n=d_n\ell_{n-1}+\ell_{n-2}$ for $n\ge 1$. 
The slope of $c_\dd$ has the continued fraction expansion \cite{Lothaire}
\begin{equation}\label{eq001}
\alpha=[0;\,d_1+1,\,d_2,\,d_3,\dots].
\end{equation}

\begin{definition}[\cite{Glen2009}]\label{de0001}
A directive sequence $\dd$ is purely periodic of period $k$ if $d_{i+k}=d_i$ for all $i\ge 1$. The notation $\dd=(d_1,\dots,d_k)^\infty$ is used with $d_1,\dots,d_k\ge 1$.
\end{definition}

\begin{remark}[\cite{HardyWright}]\label{de0004}
For purely periodic $\dd$ the slope is $\alpha=[0;\,d_1+1,\,(d_2,\dots,d_k,d_1)^\infty]$.
\end{remark}

\begin{definition}\label{de0003}
For a positive integer $d$,
\begin{equation}\label{eqq01de0003}
M(d)=\begin{pmatrix} d & 1\\ 1 & 0\end{pmatrix},\qquad M_\dd=M(d_1)M(d_2)\cdots M(d_k).
\end{equation}
\end{definition}
The matrix $M_\dd$ has determinant $(-1)^k$ and nonnegative integer entries, and some power of it is strictly positive. 
By the Perron--Frobenius theorem \cite{HornJohnson}, $M_\dd$ has a simple eigenvalue $\lambda_\dd>1$ of largest modulus. 
This growth constant $\lambda_\dd$ is a root of $x^2-\operatorname{tr}(M_\dd)\,x+(-1)^k=0$.

\subsection{Palindromic prefixes of standard Sturmian words}\label{subsec002pal}
\begin{definition}\label{def0001}
For an infinite word $x$, $L(x)$ denotes the set of lengths of the nonempty palindromic prefixes of $x$, and $\pi_x(N)=\#\{n\in L(x):\ n\le N\}$ for $N\ge 0$. For $x=c_\dd$ the notation $\pi_\dd(N)$ is used.
\end{definition}

\begin{definition}\label{def0002}
The palindromic closure $\pal{w}$ of a finite word $w$ is the shortest palindrome having $w$ as a prefix.
\end{definition}

Writing $w=uv$ with $v$ the longest palindromic suffix of $w$, the palindromic closure is $\pal{w}=uv\rev{u}$, and in particular it is unique.
For a directive sequence $\dd$ the directive word is
\begin{equation}\label{eqq0001}
\Delta_\dd=a_1a_2a_3\cdots=0^{d_1}1^{d_2}0^{d_3}1^{d_4}\cdots,
\end{equation}
in which block $i$ consists of $d_i$ equal letters and the letters alternate from block to block. 
The words $w_m$ are defined by $w_0=\varepsilon$ and $w_{m+1}=\pal{(w_ma_{m+1})}$ for $m\ge 0$.

\begin{theorem}[\cite{DeLuca,DJP}]\label{thm0closure}
The word $c_\dd$ is the limit of the words $w_m$, and the palindromic prefixes of $c_\dd$ are exactly the words $w_m$ with $m\ge 0$.
\end{theorem}

The limit statement is Theorem~5 of \cite{DeLuca}, formulated there with the palindromic left closure and the recursion $t_{n+1}=(x_nt_n)^{(-)}$ over the alphabet $\{a,b\}$. 
Droubay, Justin and Pirillo \cite{DJP} define standard episturmian words by the right closure recursion. 
In \cite{DeLuca} the directive sequence $(q_0,q_1,\dots)$ satisfies $q_{i-1}=d_i$, with $a$ read as $0$ and $b$ read as $1$. 
The length set is
\begin{equation}\label{eqq0002}
\Lambda(\dd)=\{\, j\,\ell_{i-1}+\ell_{i-2}-2 \;:\; i\ge 1,\ 1\le j\le d_i \,\}.
\end{equation}
Let $\delta_m=|w_m|$.
\begin{lemma}[\cite{JP}]\label{Lem01sec05}
Let $x=a_{m+1}$ for $m\ge 0$. 
If $x$ does not occur in $a_1\cdots a_m$, then $\delta_{m+1}=2\delta_m+1$. 
Otherwise $p=\max\{q\le m: a_q=x\}$ exists and $\delta_{m+1}=2\delta_m-\delta_{p-1}$.
\end{lemma}
Lemma~\ref{Lem01sec05} is the length form of Justin's formula $\pal{(ux)}=uu'^{-1}u$, where $u=u'xu''$ and $x$ does not occur in $u''$ \cite{Justin05}. 
It is equivalent to the relations $u_{n+1}=u_nx_nu_n$ and $u_{n+1}=u_nu_{P(n)}^{-1}u_n$ of \cite[p.~287]{JP} with $u_n=w_{n-1}$; see also \cite[Def.~4.18]{Fischler}, where $n_i=\delta_{i-1}$.

\begin{theorem}[\cite{Fischler}]\label{thm0A1}
For every directive sequence $\dd$ with $d_i\ge 1$ for all $i\ge 1$,
\begin{equation}\label{eqq01thm0A1}
L(c_\dd)=\Lambda(\dd)\setminus\{0\},
\end{equation}
where the value $0$, attained only for $i=1$ and $j=1$, is the length of the empty prefix.
\end{theorem}
A derivation from Theorem~\ref{thm0closure}, with the indexing of the lengths by the pairs $(i,j)$, is given in Section~\ref{sec01A1}.
 
\section{Main results}\label{mainsec03}

This section has two parts. Subsection~\ref{sec01A1} derives the lengths of the palindromic prefixes of $c_\dd$ from the closure construction and shows that they are indexed bijectively by the pairs $(i,j)$. Subsection~\ref{sec00Asymptotic} uses this indexing to count the palindromic prefixes of length at most $N$ and proves the asymptotic formula for purely periodic $\dd$.

\subsection{The length set of the palindromic prefixes}\label{sec01A1}

Throughout this subsection $d_i\ge 1$ for all $i\ge 1$, which holds automatically for every purely periodic $\dd$ by Definition~\ref{de0001}. In the length set $\Lambda(\dd)$ of~\eqref{eqq0002} we write $\theta_{i,j}:=j\,\ell_{i-1}+\ell_{i-2}-2$.
The symbols $D_0=0$ and $D_i=d_1+\cdots+d_i$ mark the block boundaries of the directive word \eqref{eqq0001}, in which block $i$ is the run of $d_i$ equal letters at the positions $D_{i-1}+1,\dots,D_i$. 
The words $w_m$ are given by $w_0=\varepsilon$ and $w_{m+1}=\pal{(w_ma_{m+1})}$, and $\delta_m=|w_m|$.\par 
Since $w_m$ is a proper prefix of $w_{m+1}$, the lengths satisfy $\delta_m<\delta_{m+1}$, so the $w_m$ are pairwise distinct and listed in increasing length. By Theorem~\ref{thm0closure},
\begin{equation}\label{eqq01sec5}
L(c_\dd)=\{\delta_m: m\ge 1\}.    
\end{equation}
\begin{lemma}\label{Lem02sec05}
\begin{enumerate}
\item[(i)] The map $(i,j)\mapsto m=D_{i-1}+j-1$ is a bijection from $\{(i,j): i\ge 1,\ 1\le j\le d_i\}$ onto $\N=\{0,1,2,\dots\}$, and for the $m$ corresponding to $(i,j)$ the letter $a_{m+1}$ is the $j$-th letter of block $i$.
\item[(ii)] $\theta_{i,j+1}-\theta_{i,j}=\ell_{i-1}$ for all $i\ge 1$ and $j\ge 1$.
\item[(iii)] $\theta_{i,d_i+1}=\theta_{i+1,1}$ and $\theta_{i,d_i}=\ell_i-2$ for all $i\ge 1$.
\item[(iv)] $\theta_{1,1}=0$, and $\theta_{i,j}\ge 1$ for every other pair $(i,j)$.
\end{enumerate}
\end{lemma}

\begin{proof}
For fixed $i$ the parameter $j$ runs through $1,\dots,d_i$, so $m$ runs through $D_{i-1},\dots,D_i-1$, and these intervals partition $\N$. The position $m+1=D_{i-1}+j$ is the $j$-th position of block $i$. This proves item~(i). Item~(ii) is immediate from the definition of $\theta_{i,j}$. Item~(iii) follows from $\ell_i=d_i\ell_{i-1}+\ell_{i-2}$, which gives
\begin{equation}\label{eqq01Lem02sec05}
\theta_{i,d_i}=d_i\ell_{i-1}+\ell_{i-2}-2=\ell_i-2,\qquad
\theta_{i,d_i+1}=\ell_i+\ell_{i-1}-2=\theta_{i+1,1}.
\end{equation}
For item~(iv), $\theta_{1,j}=j\ell_0+\ell_{-1}-2=j-1$, which is $0$ only for $j=1$.
Since $\ell_i\ge\ell_{i-1}+\ell_{i-2}$, the sequence $(\ell_i)_{i\ge -1}$ is nondecreasing with $\ell_{-1}=\ell_0=1$, so for $i\ge 2$ we have $\ell_{i-1}\ge\ell_1=d_1+1\ge 2$ and $\ell_{i-2}\ge 1$. Hence $\theta_{i,j}\ge 2+1-2=1$.
\end{proof}

\begin{proposition}\label{pro01sec05}
For every $i\ge 1$ and $1\le j\le d_i$, the index $m=D_{i-1}+j-1$ satisfies $\delta_m=\theta_{i,j}$.
\end{proposition}
\begin{proof}
We argue by induction on $m\ge 0$. By Lemma~\ref{Lem02sec05}, item~(i), every $m\ge 0$ corresponds to exactly one pair $(i,j)$. For $m=0$ this pair is $(1,1)$, and $\delta_0=0=\theta_{1,1}$. For the induction step the index $m\ge 0$ corresponds to $(i,j)$, the claim holds for all indices up to $m$, and $x=a_{m+1}$. It suffices to show $\delta_{m+1}=\theta_{i,j+1}$, because for $j<d_i$ the pair of $m+1$ is $(i,j+1)$, and for $j=d_i$ the value $\theta_{i,d_i+1}$ equals $\theta_{i+1,1}$, which is the value required for the pair $(i+1,1)$ of $m+1$, by Lemma~\ref{Lem02sec05}, item~(iii).

If $j\ge 2$, the letter $a_m$ lies in the same block as $x$, so $a_m=x$, the index $p$ of Lemma~\ref{Lem01sec05} equals $m$, and $\delta_{m+1}=2\delta_m-\delta_{m-1}$. The index $m-1$ corresponds to $(i,j-1)$, so the induction hypothesis and Lemma~\ref{Lem02sec05}, item~(ii) give
\begin{equation}\label{eqq01pro01sec05}
\delta_{m+1}=2\theta_{i,j}-\theta_{i,j-1}=\theta_{i,j}+\ell_{i-1}=\theta_{i,j+1}.
\end{equation}

If $j=1$ and $i=1$, then $m=0$ and $x=0$ does not occur in the empty word, so $\delta_1=2\delta_0+1=1$ by Lemma~\ref{Lem01sec05}, which equals $\theta_{1,2}=2\ell_0+\ell_{-1}-2$.

If $j=1$ and $i=2$, then $m=D_1=d_1$, the word $a_1\cdots a_m=0^{d_1}$ contains no letter $1$, and $x=1$. Lemma~\ref{Lem01sec05} and the induction hypothesis give $\delta_{m+1}=2\theta_{2,1}+1$, and $\theta_{2,1}=\ell_1+\ell_0-2=d_1$ gives $\delta_{m+1}=2d_1+1=2\ell_1+\ell_0-2=\theta_{2,2}$.

If $j=1$ and $i\ge 3$, then $m=D_{i-1}$, and $x$ is the letter of block $i$, which is also the letter of block $i-2$ and differs from the letter of block $i-1$. The last occurrence of $x$ in $a_1\cdots a_m$ is therefore the last position of block $i-2$, so $p=D_{i-2}$ in Lemma~\ref{Lem01sec05}. Then $p-1=D_{i-3}+d_{i-2}-1$ corresponds to the pair $(i-2,d_{i-2})$, and $p-1<m$, so the induction hypothesis and Lemma~\ref{Lem02sec05}, item~(iii) give $\delta_{p-1}=\theta_{i-2,d_{i-2}}=\ell_{i-2}-2$. Thus,
\begin{equation}\label{eqq02pro01sec05}
\delta_{m+1}=2\theta_{i,1}-(\ell_{i-2}-2)=2(\ell_{i-1}+\ell_{i-2}-2)-\ell_{i-2}+2=2\ell_{i-1}+\ell_{i-2}-2=\theta_{i,2}.
\end{equation}
The four situations cover every pair $(i,j)$, which completes the induction.
\end{proof}

\begin{corollary}\label{cor01sec05}
The map $(i,j)\mapsto\theta_{i,j}$ is injective on $\{(i,j): i\ge 1,\ 1\le j\le d_i\}$, and its image is $\Lambda(\dd)$. Consequently each length in $L(c_\dd)$ is attained by exactly one pair $(i,j)\ne(1,1)$, and $L(c_\dd)=\Lambda(\dd)\setminus\{0\}$.
\end{corollary}
\begin{proof}
By Lemma~\ref{Lem02sec05}, item~(i), the pairs correspond bijectively to the indices $m\ge 0$, and by Proposition~\ref{pro01sec05} the pair of $m$ has $\theta_{i,j}=\delta_m$. Since $\delta_m<\delta_{m+1}$ for all $m\ge 0$, the values $\delta_m$ are pairwise distinct, so the map is injective. By Lemma~\ref{Lem02sec05}, item~(iv), the value $0$ occurs only for $(1,1)$, and \eqref{eqq01sec5} gives the last statement.
\end{proof}
\subsection{Asymptotic count of palindromic prefixes}\label{sec00Asymptotic}
 
The sequence $\dd=(d_1,\dots,d_k)^\infty$ is purely periodic of period $k$, so $d_i\ge 1$ and $d_{i+k}=d_i$ for all $i\ge 1$. The symbol $S$ stands for $d_1+\cdots+d_k$. The only facts about palindromic prefixes used below are~\eqref{eqq01sec5}, Proposition~\ref{pro01sec05} and Corollary~\ref{cor01sec05}. By~\eqref{eqq01sec5}, $\pi_\dd(N)=\#\{m\ge 1:\ \delta_m\le N\}$. By Lemma~\ref{Lem02sec05}\,(i), the map $(i,j)\mapsto m=D_{i-1}+j-1$ is a bijection from $\{(i,j):\ i\ge 1,\ 1\le j\le d_i\}$ onto $\N$, under which $(1,1)$ corresponds to $m=0$, and by Proposition~\ref{pro01sec05}, $\delta_m=\theta_{i,j}$ for the corresponding pair. By Corollary~\ref{cor01sec05}, the pairs and the lengths correspond one to one. Therefore
\[
\pi_\dd(N)=\#\{(i,j):\ i\ge 1,\ 1\le j\le d_i,\ (i,j)\ne(1,1),\ \theta_{i,j}\le N\},
\qquad \theta_{i,j}=j\,\ell_{i-1}+\ell_{i-2}-2 .
\]
 
\begin{theorem}\label{Thm01Asymptotic}
For every purely periodic $\dd$ of period $k$,
\begin{equation}\label{eqq01Thm01Asymptotic}
\pi_\dd(N)\sim\frac{d_1+\cdots+d_k}{\log\lambda_\dd}\,\log N .
\end{equation}
\end{theorem}
\begin{proof}
The vectors $v_n=(\ell_n,\ell_{n-1})^T$ satisfy $v_0=(1,1)^T$ and $v_n=M(d_n)v_{n-1}$ for $n\ge 1$. Since $d_{pk+s}=d_s$ for $1\le s\le k$, the product $M(d_{qk})M(d_{qk-1})\cdots M(d_1)$ is the $q$-th power of $B=M(d_k)\cdots M(d_1)$, and therefore $v_{qk}=B^qv_0$. Each $M(d)$ is symmetric with determinant $-1$, so $B$ is the transpose of $M_\dd$, with $\det B=(-1)^k$ and $\operatorname{tr}B=t$, where $t$ denotes $\operatorname{tr}M_\dd$. The Cayley--Hamilton theorem gives $B^2=tB-(-1)^kI$. Multiplying by $B^qv_0$ and reading the first coordinate gives
\begin{equation}\label{eqq02Thm01Asymptotic}
y_{q+2}=t\,y_{q+1}-(-1)^k y_q,
\end{equation}
for $q\ge 0$, where $y_q:=\ell_{qk}$, $y_0=1$ and $y_1=\ell_k$.
 
\noindent The matrices $A_n=M(d_1)\cdots M(d_n)$, with $A_0=I$, are written
\[
A_n=\begin{pmatrix}P_n&P_{n-1}\\ Q_n&Q_{n-1}\end{pmatrix},
\qquad P_{-1}=0,\ P_0=1,\ Q_{-1}=1,\ Q_0=0 .
\]
The relation $A_n=A_{n-1}M(d_n)$ gives, for all $n\ge 1$,
\begin{equation}\label{eqq03Thm01Asymptotic}
P_n=d_nP_{n-1}+P_{n-2}, \qquad Q_n=d_nQ_{n-1}+Q_{n-2}.
\end{equation}
Since $P_1=d_1\ge 1$ and $Q_1=1$, induction on $n$ gives $P_n\ge 1$ for $n\ge 0$ and $Q_n\ge 1$ for $n\ge 1$. Since $A_k=M_\dd$, the trace is $t=P_k+Q_{k-1}$. We claim that
\begin{equation}\label{eqq03bThm01Asymptotic}
t\ge 1\ \text{ for every } k\ge 1,\qquad t\ge 3\ \text{ if } k\ge 2 .
\end{equation}
For $k=1$ we have $t=P_1+Q_0=d_1\ge 1$. For $k\ge 2$, \eqref{eqq03Thm01Asymptotic} gives $P_k=d_kP_{k-1}+P_{k-2}\ge P_{k-1}+P_{k-2}\ge 2$ and $Q_{k-1}\ge 1$, hence $t\ge 3$.
 
\noindent The polynomial $f(x)=x^2-tx+(-1)^k$ satisfies $f(1)=1-t+(-1)^k$. If $k$ is odd, then $f(1)=-t<0$ by the first part of~\eqref{eqq03bThm01Asymptotic}. If $k$ is even, then $k\ge 2$, and $f(1)=2-t<0$ by the second part of~\eqref{eqq03bThm01Asymptotic}. In both cases $f(1)<0$.
Consequently $f$ has two real roots with $1$ strictly between them. The larger root $\lambda$ exceeds $1$, and the other root is $\mu=(-1)^k/\lambda$, $|\mu|=1/\lambda<1$ and $\lambda\ne\mu$. Both roots are eigenvalues of $M_\dd$, since $f$ is its characteristic polynomial, and $\lambda$ has the larger modulus. Thus, $\lambda=\lambda_\dd$.
 
\noindent The sequence $\ell_n$ is strictly increasing for $n\ge 0$: for $n\ge 1$ it holds that $\ell_n=d_n\ell_{n-1}+\ell_{n-2}\ge\ell_{n-1}+\ell_{n-2}>\ell_{n-1}$, where $\ell_{n-2}\ge 1$ and $\ell_{-1}=\ell_0=1$. Since $\lambda\ne\mu$, there are constants $\kappa_1,\kappa_2$ with $\kappa_1+\kappa_2=1$ and $\kappa_1\lambda+\kappa_2\mu=\ell_k$. The sequence $\kappa_1\lambda^q+\kappa_2\mu^q$ has the same first two terms as $y_q$ and satisfies the same recurrence, because $\lambda^2=t\lambda-(-1)^k$ and $\mu^2=t\mu-(-1)^k$. Therefore
\begin{equation}\label{eqq03cThm01Asymptotic}
y_q=\ell_{qk}=\kappa_1\lambda^q+\kappa_2\mu^q\qquad\text{for all } q\ge 0 .
\end{equation}
Subtracting $\mu$ times the equation $\kappa_1+\kappa_2=1$ from $\kappa_1\lambda+\kappa_2\mu=\ell_k$ gives $\kappa_1(\lambda-\mu)=\ell_k-\mu$. Since $\lambda-\mu>0$ and $\ell_k\ge\ell_1=d_1+1\ge 2>|\mu|$, we get $\kappa_1>0$. Thus, from~\eqref{eqq03cThm01Asymptotic},
 
\begin{equation}\label{eqq04Thm01Asymptotic}
\ell_{qk}=\kappa_1\lambda^q\Bigl(1+\tfrac{\kappa_2}{\kappa_1}\bigl(\tfrac{\mu}{\lambda}\bigr)^q\Bigr),
\qquad\text{and}\qquad \frac{\log\ell_{qk}}{q}\longrightarrow\log\lambda\quad(q\to\infty),
\end{equation}
where $|\mu/\lambda|<1$.
 
\noindent For $n\ge k$ the integer $q=\lfloor n/k\rfloor$ satisfies $qk\le n<(q+1)k$, and monotonicity of $\ell_n$ gives $\ell_{qk}\le\ell_n\le\ell_{(q+1)k}$. Therefore
\begin{equation}\label{eqq05Thm01Asymptotic}
\frac{\log\ell_{qk}}{q}\cdot\frac{q}{(q+1)k}\ \le\ \frac{\log\ell_n}{n}\ \le\ \frac{\log\ell_{(q+1)k}}{q+1}\cdot\frac{q+1}{qk},
\end{equation}
and both bounds tend to $(\log\lambda)/k$ as $n\to\infty$. Thus, \eqref{eqq04Thm01Asymptotic} and \eqref{eqq05Thm01Asymptotic} give
\begin{equation}\label{eqq06Thm01Asymptotic}
\lim_{n\to\infty}\frac{\log \ell_n}{n}=\frac{\log\lambda_\dd}{k}.
\end{equation}
 
\noindent For $N\ge 0$ the index $I=I(N)=\max\{i\ge 0:\ \ell_i\le N+2\}$ exists, since $\ell_0=1$ and $\ell_i\to\infty$, and $I(N)\to\infty$ as $N\to\infty$. For $N$ large enough that $I\ge 1$, the definition of $I$ gives $\ell_I\le N+2<\ell_{I+1}$,
\begin{equation}\label{eqq07Thm01Asymptotic}
\frac{\log\ell_I}{I}\ \le\ \frac{\log(N+2)}{I}\ <\ \frac{\log\ell_{I+1}}{I+1}\cdot\frac{I+1}{I}.
\end{equation}
Both outer terms tend to $(\log\lambda)/k$. Thus, $I(N)/\log(N+2)\to k/\log\lambda$, and since $\log(N+2)/\log N\to 1$,
\[
\frac{I(N)}{\log N}\longrightarrow\frac{k}{\log\lambda}.
\]
 
\noindent The set of pairs counted by $\pi_\dd(N)$ is bounded from both sides as follows. For $1\le i\le I$ and $1\le j\le d_i$, the identity $\theta_{i,j+1}-\theta_{i,j}=\ell_{i-1}>0$ gives $\theta_{i,j}\le\theta_{i,d_i}=\ell_i-2\le\ell_I-2\le N$. Hence every pair with $i\le I$ is counted, except $(1,1)$, and
\begin{equation}\label{eqq08Thm01Asymptotic}
\pi_\dd(N)\ge\sum_{i=1}^{I}d_i-1.
\end{equation}
In the other direction, let $(i,j)$ be a pair with $\theta_{i,j}\le N$. Since $\theta_{i,1}\le\theta_{i,j}$, we get $\ell_{i-1}+\ell_{i-2}-2\le N$, hence $\ell_{i-1}\le N+2$, and so $i-1\le I$ by the definition of $I$. Thus every counted pair has $i\le I+1$, and
\begin{equation}\label{eqq08bThm01Asymptotic}
\pi_\dd(N)\le\sum_{i=1}^{I+1}d_i .
\end{equation}
 
\noindent Periodicity of $\dd$ gives $\sum_{i=1}^{m}d_i=S\lfloor m/k\rfloor+r_m$ with $0\le r_m\le S$, so $\bigl|\sum_{i=1}^{m}d_i-Sm/k\bigr|\le S$ for every $m\ge 0$. Together with~\eqref{eqq08Thm01Asymptotic} and~\eqref{eqq08bThm01Asymptotic},
\begin{equation}\label{eqq09Thm01Asymptotic}
\frac{S\,I}{k}-S-1\ \le\ \pi_\dd(N)\ \le\ \frac{S\,(I+1)}{k}+S .
\end{equation}
 
\noindent Divide~\eqref{eqq09Thm01Asymptotic} by $\log N$ and use $I(N)/\log N\to k/\log\lambda_\dd$. Both outer terms then tend to $S/\log\lambda_\dd$, since $\frac{S}{k}\cdot\frac{k}{\log\lambda_\dd}=\frac{S}{\log\lambda_\dd}$. Hence
\begin{equation}\label{eqq10Thm01Asymptotic}
\frac{\pi_\dd(N)}{\log N}\longrightarrow\frac{S}{\log\lambda_\dd},
\end{equation}
which is~\eqref{eqq01Thm01Asymptotic}.
\end{proof}

\begin{corollary}\label{cornu01}
For every purely periodic $\dd$ of period $k$,
\[
\pi_\dd(N)=\frac{S}{\log\lambda_\dd}\,\log N+O(1).
\]
\end{corollary}
\begin{proof}
By~\eqref{eqq03cThm01Asymptotic} and $\kappa_1>0$, there are constants $c_1,c_2>0$ with $c_1\lambda^q\le\ell_{qk}\le c_2\lambda^q$ for all $q\ge 0$, where $\lambda=\lambda_\dd$. For $n\ge 0$ and $q=\lfloor n/k\rfloor$, monotonicity of $\ell_n$ gives $\ell_{qk}\le\ell_n\le\ell_{(q+1)k}$, and $\lambda^{n/k-1}<\lambda^q\le\lambda^{n/k}$. Hence $C_1\lambda^{n/k}\le\ell_n\le C_2\lambda^{n/k}$ with $C_1=c_1/\lambda$ and $C_2=c_2\lambda$.
From $\ell_I\le N+2<\ell_{I+1}$ we obtain
\[
C_1\lambda^{I/k}\le N+2,\qquad N+2<C_2\lambda^{(I+1)/k},
\]
that is,
\[
\frac{k\log\bigl((N+2)/C_2\bigr)}{\log\lambda}-1< I\le\frac{k\log\bigl((N+2)/C_1\bigr)}{\log\lambda}.
\]
Since $\log(N+2)-\log N\le\log 3$ for $N\ge 1$, there is a constant $C$ with $\bigl|I(N)-k\log N/\log\lambda\bigr|\le C$ for all $N\ge 1$. Substituting this into~\eqref{eqq09Thm01Asymptotic} gives
\[
\frac{S}{\log\lambda}\log N-\frac{SC}{k}-S-1\le\pi_\dd(N)\le\frac{S}{\log\lambda}\log N+\frac{SC}{k}+\frac{S}{k}+S,
\]
which proves the claim.
\end{proof}
 
\section{Discussion}

The description of the lengths of palindromic prefixes in Theorem~\ref{thm0A1} agrees with the one stated in Section~3.1 of Fischler's study of palindromic prefixes \cite{Fischler}, where it is attributed to de Luca \cite{DeLuca} and to Droubay, Justin and Pirillo \cite{DJP}. Section~\ref{sec01A1} derives it directly from the closure recursion and adds the explicit bijection $(i,j)\mapsto m$ of Lemma~\ref{Lem02sec05} and Corollary~\ref{cor01sec05}, which is what makes the counting possible. The content specific to this paper is the counting statement. We are not aware of an earlier statement of the asymptotic formula for $\pi_\dd(N)$.

In the literature the abundance of palindromic prefixes is measured by $\limsup n_{i+1}/n_i$, where $n_i$ is the length of the $i$-th palindromic prefix \cite{Fischler}. By Proposition~\ref{pro01sec05}, the ratios of consecutive lengths are governed by the quotients $\ell_{i-2}/\ell_{i-1}=[0;d_{i-1},\dots,d_2,d_1+1]$, so this quantity records the largest relative gap between consecutive palindromic prefixes and depends on the extreme values of the directive sequence. The count $\pi_\dd(N)$ takes all gaps into account, and its constant $S/\log\lambda_\dd$ depends on $\dd$ in a different way.

This constant is the ratio of $S/k$, the mean of $d_1,\dots,d_k$, to $(\log\lambda_\dd)/k$, the exponential growth rate of $\ell_n$ per index, by~\eqref{eqq06Thm01Asymptotic}. It is not a function of the mean alone, because $\lambda_\dd$ depends on the order of $d_1,\dots,d_k$ and not only on their sum. For $\dd=(2,2)^\infty$ we have $\lambda_\dd=(1+\sqrt2)^2$ and the constant is $2/\log(1+\sqrt2)\approx 2.27$, while for $\dd=(1,3)^\infty$ we have $\lambda_\dd=(5+\sqrt{21})/2$ and the constant is $4/\log\bigl((5+\sqrt{21})/2\bigr)\approx 2.55$, although both sequences have mean $2$. Replacing $k$ by $mk$ multiplies $S$ by $m$ and replaces $M_\dd$ by $M_\dd^{\,m}$, so $\log\lambda_\dd$ is also multiplied by $m$, and the constant does not depend on the choice of period. For $\dd=(1)^\infty$ the theorem gives $\pi_\dd(N)\sim\log N/\log\varphi$ with $\varphi=(1+\sqrt5)/2$.

Periodicity enters the proof of Theorem~\ref{Thm01Asymptotic} only through~\eqref{eqq06Thm01Asymptotic} and the bounded oscillation of the partial sums of $d_i$. The remaining steps hold for every directive sequence, which gives the following statement.

\begin{proposition}\label{prop01disc}
Let $\dd$ be a directive sequence with $d_i\ge 1$ for all $i\ge 1$, and let $I(N)=\max\{i\ge 0:\ \ell_i\le N+2\}$. Then
\[
\sum_{i\le I(N)}d_i-1\le\pi_\dd(N)\le\sum_{i\le I(N)+1}d_i .
\]
If $\frac{\log\ell_n}{n}\to\rho>0$ and $\frac1n\sum_{i\le n}d_i\to\bar d\in(0,\infty)$, then $\pi_\dd(N)\sim(\bar d/\rho)\log N$.
\end{proposition}
\begin{proof}
The two bounds are~\eqref{eqq08Thm01Asymptotic} and~\eqref{eqq08bThm01Asymptotic}, whose proofs use only Lemma~\ref{Lem02sec05} and the monotonicity of $\ell_n$. By~\eqref{eqq07Thm01Asymptotic} and the assumption on $\log\ell_n/n$, we get $I(N)/\log N\to 1/\rho$. Since $\sum_{i\le I}d_i=(\bar d+o(1))I$, dividing the bounds by $\log N$ gives the claim.
\end{proof}

The error term in Corollary~\ref{cornu01} cannot be improved to $o(1)$. Indeed, the lengths $\theta_{i,j}$ are pairwise distinct by Corollary~\ref{cor01sec05}, so $\pi_\dd(N)$ increases by exactly $1$ at each of them, while $\frac{S}{\log\lambda_\dd}\log N$ changes by $o(1)$ across such a point. Hence the difference $\pi_\dd(N)-\frac{S}{\log\lambda_\dd}\log N$ has no limit.

Finally, the periodic case is special. For $n\ge 1$ the number $\ell_n$ is the denominator of the $n$-th convergent of $\alpha=[0;d_1+1,d_2,\dots]$, so by L\'evy's theorem we have $\log\ell_n/n\to\rho=\pi^2/(12\log 2)$ for almost every $\alpha\in (0,1/2)$ \cite{IK02,DK02}. By the Borel--Bernstein theorem \cite{IK02}, for almost every $\alpha$ the inequality $d_n\ge n\log n$ holds for infinitely many $n$. For such $n$ and $N=\ell_n-2$ we have $I(N)=n$, so Proposition~\ref{prop01disc} gives $\pi_\dd(N)\ge d_n-1\ge n\log n-1$, while $\log N<\log\ell_n\sim\rho n$. Hence $\pi_\dd(N)/\log N$ is unbounded, and no asymptotic of the form $C\log N$ holds for such $\dd$.

\section{Conclusion}

We derived the lengths of the palindromic prefixes of a standard Sturmian word $c_\dd$ with directive sequence $\dd$, where $d_i\ge 1$ for all $i\ge 1$, directly from the construction by iterated palindromic closure. They are the numbers $j\,\ell_{i-1}+\ell_{i-2}-2$ with $i\ge 1$ and $1\le j\le d_i$, where $\ell_n$ is the length of the $n$-th standard word. The value $0$ occurs only for $(i,j)=(1,1)$ and corresponds to the empty prefix, and every nonzero length is attained by exactly one pair $(i,j)$. The proof tracks the lengths of the closures through the block boundaries of the directive word $\Delta_\dd$.

For a purely periodic directive sequence $\dd=(d_1,\dots,d_k)^\infty$, the counting function $\pi_\dd(N)$ of nonempty palindromic prefixes of length at most $N$ satisfies
\[
\pi_\dd(N)\sim\frac{d_1+\cdots+d_k}{\log\lambda_\dd}\,\log N,
\]
where $\lambda_\dd$ is the largest eigenvalue of $M(d_1)\cdots M(d_k)$. More precisely, Corollary~\ref{cornu01} gives $\pi_\dd(N)=\frac{d_1+\cdots+d_k}{\log\lambda_\dd}\log N+O(1)$, and this error term cannot be improved to $o(1)$. The constant does not depend on the choice of period, and for the Fibonacci word it equals $1/\log\varphi$. The argument rests on two ingredients, the exponential growth of $\ell_n$ governed by $\lambda_\dd$ and the bounded oscillation of the partial sums of $d_i$ over one period. The constant is determined by the sum of the $d_i$ over one period together with the growth rate of $\ell_n$, and not by the mean of the $d_i$ alone. It differs in nature from the quantity $\limsup n_{i+1}/n_i$, which is governed by the extreme values of the quotients $\ell_{i-2}/\ell_{i-1}$.

Proposition~\ref{prop01disc} extends the counting argument to arbitrary directive sequences with $d_i\ge 1$, and together with the Borel--Bernstein theorem shows that $\pi_\dd(N)/\log N$ is unbounded for almost every slope.
 

\end{document}